\documentclass[a4paper,11pt]{amsart}

\usepackage{amsfonts}
\usepackage{latexsym}
 \usepackage{amssymb}
\usepackage[all]{xy}

\newtheorem{thm}{\sc theorem}[section]

\newtheorem{lemma}[thm]{\sc lemma}

\newenvironment{pf}{{{\bf Proof.}}}{\hfill$\Box$\\[1mm]}

\newcommand{\sm}{{\sigma_\mathrm{mod}}}

\renewcommand{\b}{\mathsf{b}} 
\renewcommand{\d}{\mathsf{d}}

\renewcommand{\k}{k}

\renewcommand{\t}{\mathsf{t}}

\begin{document}

\title{The Hochschild cohomology ring of the standard 
Podle\'s quantum sphere}

\author{Ulrich Kr\"ahmer} 

\address{University Gardens, Glasgow G12 8QW, UK}

\email{ukraehmer@maths.gla.ac.uk}

\begin{abstract}
The cup and cap product in twisted 
Hochschild
(co)\-ho\-mo\-lo\-gy is computed 
for the standard quantum
2-sphere and used 
to construct a cyclic
2-cocycle that represents the 
fundamental Hochschild class.
\end{abstract}
\maketitle
 
\section{Introduction}
The aim of this 
article is to compute for the quantised
coordinate ring $A=\mathbb{C}_q[S^2]$ 
(the standard Podle\'s quantum 2-sphere) the
cup and cap product 
\begin{eqnarray}
		  \smallsmile 
&:& H^m(A,{}_\sigma A) \otimes 
		  H^n(A,{}_\tau A) \rightarrow 
		  H^{m+n}(A,{}_{\tau \circ \sigma} A),\nonumber\\ 
		  \smallfrown 
&:& H_n(A,{}_\sigma A) \otimes 
		  H^m(A,{}_\tau A) \rightarrow
		  H_{n-m} (A,{}_{\tau \circ \sigma}
		  A)\nonumber
\end{eqnarray} 
in the Hochschild (co)homology
of $A$ with coefficients in the bimodules 
${}_\sigma A$ that arise
by twisting the canonical bimodule 
structure of $A$
by $ \sigma \in \mathrm{Aut} (A)$. 
In \cite{tomuli2} we carried out similar 
computations for $\mathbb{C}_q[SU(2)]$. 
For $\mathbb{C}_q[S^2]$ 
they become much simpler
and their conceptual content that was
somewhat hidden in 
\cite{tomuli2} between lengthy 
computations becomes more transparent.

For
the coordinate ring of a smooth variety 
$X$, the 
Hochschild-Kostant-Rosenberg 
theorem identifies 
$\Lambda(A):=
(\bigoplus_{n \ge 0} H^n(A,A),\smallsmile)$ 
with the algebra of multivector fields on $X$, 
and $\bigoplus_{n \ge 0} H_n(A,A)$
as a $\Lambda(A)$-module (via $\smallfrown$)  
with the differential forms $ \Omega (X)$ on $X$. 
For noncommutative algebras, 
$\Lambda(A)$ tends to be 
fairly degenerate. However, twisting by 
$ \sigma $
allows one to consider richer cohomology 
rings that encode more
information about $A$.

For example, the
generators of the Drinfeld-Jimbo quantisation
of the Lie algebra of $SU(2)$
act via twisted rather than usual derivations 
on $A=\mathbb{C}_q[S^2]$. 
These give rise to 
two cohomology classes 
$[\partial_{\pm 1}] \in 
H^1(A,{}_{\sigma^{-1}_\mathrm{mod}} A)$, 
where $\sm$ is Woronowicz's 
modular automorphism determined 
for example by (\ref{podushka}) below. 
We will see that 
they behave under the cup product 
similar to the corresponding 
classical $SU(2)$-invariant
vector fields on 
$S^2=SU(2)/U(1)$,
$$
		  [\partial_1] \smallsmile [\partial_1]=
		  [\partial_1] \smallsmile [\partial_{-1}]+
		  q^2 [\partial_{-1}] \smallsmile [\partial_1]=
		  [\partial_{-1}] \smallsmile [\partial_{-1}]=0,
$$
and use them to define
a functional 
on $H_2(A,{}_\sm A) \simeq \mathbb{C}$ 
of the form
\begin{equation}\label{blinde}
		  \varphi ([\omega]):= 
		  q^{-1} \int [\omega] \smallfrown
  		  ([\partial_1] \smallsmile [\partial_{-1}]) 
		  \in \mathbb{C}.
\end{equation}  
Here $[\omega] \in H_2(A,{}_\sm A)$ 
is acted on by 
$[\partial_1] \smallsmile [\partial_{-1}]
\in H^2(A,{}_{\sigma^{-2}_\mathrm{mod}} A)$
to produce a class in 
$H_0(A,{}_{\sigma^{-1}_\mathrm{mod}} A)$, 
and then one applies a certain twisted
trace $\int \in (H_0(A,{}_\sigma A))^*$
in order to obtain a numerical invariant of
$[\omega]$.

The functional $\varphi$ 
provides a dual 
description of the fundamental class 
$[\d A] \in H_2(A,{}_\sm A)$ that corresponds 
under the Poincar\'e-type duality \cite{ulineu}
\begin{equation}\label{podushka}
		  H_n(A,{}_\sm A) \simeq H^{2-n}(A,A)
\end{equation} 
to 
$1 \in H^0(A,A)$. 
From the practical point of view, one can
use $\varphi$ to determine the homology class 
of a given Hochschild cycle, and 
for $\mathbb{C}_q[SU(2)]$ this tool 
allowed us
to compute the 
cyclic homology \cite{kmt} 
built upon $H_n(A,{}_\sigma A)$
as a special case of Connes-Moscovici's 
Hopf-cyclic homology \cite{cm}.

The trace $\int$ in (\ref{blinde}) is for 
$\mathbb{C}_q[S^2]$ actually a character
(namely the restriction of the counit 
$\varepsilon$ of $\mathbb{C}_q[SU(2)]$ 
to $\mathbb{C}_q[S^2]$), so 
on the level of 
chains $a_0 \otimes a_1 \otimes a_2$
in the standard Hochschild complex, 
$\varphi$ acts as
$$
		  \varphi (a_0,a_1,a_2)=
		  q^{-1}\varepsilon(a_0)F(a_1)E(a_2),
$$
where $E,F : A \rightarrow k$ are 
the (untwisted) derivations given by
$$
		  E(a):=\varepsilon(\partial_{-1}(a)),
		  \quad
		  F(a):=\varepsilon(\partial_1(a)).
$$
There seems to be a general 
principle behind this that we observed already  
in \cite{tomuli2}. Therein, the trace $\int$
needed was the integral over 
the unquantised maximal torus in quantum 
$SU(2)$. For the quantum 2-sphere 
the corresponding submanifold of $S^2$
is simply a point, namely one of the two
leaves of the symplectic foliation of the 
Poisson manifold $S^2$ quantised by $A$.

Finally we discuss how to add a 
counter term $\eta$ to $\varphi$ in order 
to obtain a functional on cyclic 
homology without changing   
the functional on $H_2(A,{}_\sm A)$.
Schm\"udgen and Wagner have defined a 
cyclic 2-cocycle in \cite{sw} that looks
like $ \varphi $, only the trace 
$\int$ is the Haar functional of 
$\mathbb{C}_q[SU(2)]$, and this makes
their functional trivial on
Hochschild homology \cite{tom}.

The structure of the paper is as follows.
In Sections~2-5 we recall the definition
of Hochschild (co)homology, of the cup and 
cap product, give some more 
background about the above mentioned 
Poincar\'e duality and introduce then 
the algebra $A=\mathbb{C}_q[S^2]$
we want to study.
In Section~6 we recall from \cite{tom}
an explicit formula for the fundamental 
class $\d A$ of $A$. In Section~7 
we determine the twisted
centre of $A$ and its cap 
product action on the second 
Hochschild homology: in \cite{tom}
it was shown that 
$H_2(A,{}_\sigma A)=0$ except when 
$\sigma=\sigma^n_\mathrm{mod}$ 
for some $n \ge 1$, and then 
$H_2(A,{}_{\sigma^n_\mathrm{mod}} A) \simeq
\mathbb{C}$. Here 
we identify the sum of all these
nontrivial homology groups with the free
module of rank 1 over a polynomial ring 
$k[x_0]$ that constitutes the twisted centre of 
$A$. We continue in Section~8 with  
recalling from \cite{tom} 
(but in a slightly simplified form) 
the computation of 
the zeroth Hochschild homology groups of 
$A$ and of the twisted 
traces on $A$, describing in addition
the cap product action of the twisted centre.
Section~9 is the first really interesting one, 
here we compute the cup product between 
twisted derivations of $A$ that arise from the 
action of the quantised Lie algebra of 
$SU(2)$ on $A$. We observe 
similar as in \cite{tomuli2} that these 
twisted derivations generate a quantised 
exterior algebra. These computations 
are then used in Section~10 to 
define and discuss $ \varphi $. 
Section~11 recalls the definition and some key 
properties of cyclic homology, and finally 
we show in Section~12 that $ \varphi $ can be 
altered by a Hochschild 
coboundary to obtain a cyclic cocycle.\\    

I acknowledge support by the EPSRC 
fellowship EP/E/043267/1. 

\section{Hochschild (co)homology with coefficients in
 ${}_\sigma A$}\label{twihoco}
Let $A$ be a unital associative algebra over a field $\k$
and $ \sigma \in \mathrm{Aut} (A)$ be an 
automorphism. 
We denote by ${}_\sigma A$ the $A$-bimodule 
which is $A$ as vector space with left and right
$A$-actions given by
$a \triangleright b \triangleleft c:=\sigma(a)bc$,
$a,b,c \in A$, and by $H_n(A,{}_\sigma A)$
and $H^n(A,{}_\sigma A)$ the Hochschild
(co)homology groups of $A$ with coefficients in 
${}_\sigma A$. Explicitly, 
$H_n(A,{}_\sigma A)$ is the homology
of the chain complex 
$C^\sigma_n:=A^{\otimes n +1}$ with boundary map    
$\b_n : C^\sigma_n \rightarrow C^\sigma_{n-1}$
given by
\begin{eqnarray}
		  \b_n(a_0 \otimes \ldots \otimes a_n)  
&=& \sum_{i=0}^{n-1} (-1)^i 
		  a_0 \otimes \ldots \otimes a_ia_{i+1} \otimes
		  \ldots \otimes a_n\nonumber\\ 
&& +(-1)^n \sigma (a_n) a_0 \otimes \ldots \otimes a_{n-1}.\nonumber  
\end{eqnarray} 
Usually we will write  
$\b(a_0,\ldots,a_n)$ instead of
$\b_n(a_0 \otimes \ldots \otimes a_n)$ 
and similarly for
other multilinear maps. 
Dually, $H^n(A,{}_\sigma A)$
is the cohomology of the cochain complex 
$C^n_\sigma:=\mathrm{Hom}_k(A^{\otimes n},A)$
with coboundary map given by
\begin{eqnarray}
		  (\b^n \psi)(a_0,\ldots,a_n) 
&=& \sigma (a_0) \psi (a_1,\ldots,a_n) 
		  \nonumber\\ 
&& +  \sum_{i=0}^{n-1} (-1)^{i+1} 
		  \psi(a_0,\ldots,a_ia_{i+1},\ldots,a_n) 
		  \nonumber\\ 
&& +(-1)^{n+1} \psi 
		  (a_0,\ldots,a_{n-1})a_n.\nonumber  
\end{eqnarray}

In degree 0, we identify 
$ \psi : A^{\otimes 0}:=k \rightarrow A$ with 
$ a:=\psi (1) \in A$. 
This is a cocycle precisely when 
$ab=\sigma(b)a$ for all $b \in A$. Thus 
$H^0(A,{}_\sigma A)$ consists of the 
$ \sigma $-central elements of $A$. In 
degree 1, a cocycle is a 
$ \sigma $-twisted derivation
$ \psi : A \rightarrow A$, $ \psi (ab)=\sigma (a)
\psi (b) + \psi (a)b$, 
and $H^1(A,{}_\sigma A)$
is the space of all such 
derivations modulo those of
the form $ \psi (a)=ba-\sigma(a)b$ 
for some $b \in A$.  
For more information and details, see 
e.g.~\cite{ce,tomuli,kmt,loday,weibel}.

\section{The cup and cap product}
The cup product is the map
$$
		  \smallsmile \,\,: H^m(A,{}_\sigma A) \otimes 
		  H^n(A,{}_\tau A) \rightarrow 
		  H^{m+n}(A,{}_{\tau \circ \sigma} A),\quad
		  \sigma , \tau \in \mathrm{Aut} (A)
$$ 
given on the level of cochains by
\begin{equation}\label{verz}
		  (\varphi \smallsmile \psi) 
		  (a_1,\ldots,a_{m+n})=
		  \tau (\varphi (a_1,\ldots,a_m)) 
		  \psi (a_{m+1},\ldots,a_{m+n}).
\end{equation} 
For any monoid 
$G \subset \mathrm{Aut} (A)$, it
turns 
$$
		  \Lambda_G(A):=
		  \bigoplus_{n \in \mathbb{N},\sigma \in G}
		  H^n(A,{}_\sigma A)
$$ 
into an $\mathbb{N} \times G$-graded 
algebra that we would like to view as some analogue of
an algebra of multivector fields on a classical space.
Twisted derivations play here the role of vector
fields, and the following easily checked 
(see \cite{tomuli2}) relations demonstrate  
their behaviour under $\smallsmile\,$: 
\begin{lemma}\label{relands}
In degree 0, $\smallsmile$
reduces to the opposite product of $A$,
\begin{equation}\label{oppo}
		  a \smallsmile b=ba,\quad
		  a \in H^0(A,{}_\sigma A),
		  b \in H^0(A,{}_\tau A),
\end{equation} 
and for 
$c \in H^0(A,{}_\sigma A)$ and twisted 
derivations $\varphi \in C^1_\sigma,
\psi \in C^1_\tau$
we have
$$
		  \psi \smallsmile c= \tau(c) \smallsmile
			\psi,\quad
		  [\varphi] \smallsmile [\psi]=
		  -[\sigma ^{-1} \circ \psi \circ \sigma]
		  \smallsmile [\varphi] \in 
		  H^2(A,{}_{\tau \circ \sigma} A).
$$
\end{lemma}

Dually,
$$
		  \Omega^G(A):=
		  \bigoplus_{n \in \mathbb{N},\sigma \in G}
		  H_n(A,{}_{\sigma^{-1}} A)
$$
becomes an $ \mathbb{N} \times G$-graded 
(right) module over $ \Lambda_G(A)$
via the cap product
$$
		  \smallfrown \,\,: H_n(A,{}_\sigma A) \otimes 
		  H^m(A,{}_\tau A) \rightarrow
		  H_{n-m} (A,{}_{\tau \circ \sigma} A),\quad m \le n.
$$
Explicitly, this is given between  
$ \varphi \in C^m_\tau$ and 
$a_0 \otimes \ldots \otimes a_n \in C_n^\sigma$ by
$$
		  (a_0 \otimes \ldots \otimes a_n) \smallfrown
		  \varphi =
		  \tau (a_0) \varphi (a_1,\ldots,a_m) \otimes 
		  a_{m+1} \otimes \ldots \otimes a_n \in
		  C_{n-m}^{\tau \circ \sigma}. 
$$
In particular, the cap product
with a twisted central
element 
$c \in H^0(A,{}_\sigma A)$ is simply given by
multiplication from the \emph{left},
\begin{equation}\label{spaetzle}
		  (a_0 \otimes \ldots \otimes a_n) \smallfrown c=
		  \sigma (a_0)c \otimes \ldots \otimes a_n=
		  ca_0 \otimes \ldots \otimes a_n.
\end{equation} 
 For more information, 
see e.g.~\cite{ce,tomuli2,nesttsygan}. 

\section{Poincar\'e duality}
The cup and cap 
product structures are intimately
related to Poincar\'e-type dualities
between homology and 
cohomology. As became clear in
recent years, there is for many algebras $A$ a
distinguished automorphism 
$ \sm $ and for all 
$ \tau \in \mathrm{Aut} (A)$ a canonical
$k$-linear isomorphism
\begin{equation}\label{podu}
		  H^n(A,{}_\tau A) \simeq 
		  H_{\mathrm{dim}(A)-n}(A,{}_{\tau \circ \sm} A),
\end{equation}    
where 
$ \mathrm{dim} (A)$ is the dimension of $A$ in
the sense of \cite{ce}, see 
e.g.~\cite{bz,farinati,tomuli,uli,ulineu} 
but first of all \cite{vdb} for this story. Under the above isomorphism,
the canonical element $1 \in H^0(A,A)$ corresponds  
to a class $[\d A] \in H_{\mathrm{dim}(A)}(A,{}_\sm A)$,
and then the isomorphism is given by the cap
product with this fundamental class. In \cite{tomuli2} we carried
this out explicitly for the standard quantised coordinate ring 
$\mathbb{C}_q[SU(2)]$ (see e.g.~\cite{chef} for
background on quantum groups), and the aim here is to
do the same for the standard 
quantum 2-sphere of
Podle\'s that we introduce in the next section.
For the coordinate ring of a smooth affine variety 
such a duality will hold if and only if the variety is
Calabi-Yau, that is, if the line bundle on $X$ 
whose sections are
the top degree K\"ahler
differentials $\Omega^{\mathrm{dim}(X)}(X)$ 
is trivial (in general one has to twist
not by an automorphism but by this module, see 
e.g.~\cite{uli}). This happens
if and only if there is a nowhere (i.e.~in no
localisation at prime ideals) vanishing element 
in $\Omega^{\mathrm{dim}(X)}(X)$, 
and under the
Hochschild-Kostant-Rosenberg isomorphism 
$\Omega^{\mathrm{dim}(X)}(X) \simeq
H_{\mathrm{dim}(X)}(k[X],k[X])$ such an element will be
identified with the fundamental class $[\d k[X]]$.
 
\section{The Podle\'s sphere}\label{aldlos} 
From now on we fix $\k=\mathbb{C}$, an element 
$q \in \k \!\setminus\!\!\{0\}$ assumed to be not a root
of unity, and $A$ is the standard Podle\'s quantum 2-sphere \cite{Po}, 
that is, the universal $\k$-algebra 
generated by $x_{-1},x_0,x_1$
satisfying the relations
$$
		  x_{\pm 1}x_0=q^{\mp 2}x_0x_{\pm 1},\quad
		  x_{\pm 1}x_{\mp 1}=
		  q^{\mp 2}x_0^2+q^{\mp 1}x_0.
$$
It follows easily from these relations that 
the elements 
$$
		  e_{ij}:=\left\{
\begin{array}{ll}
x_0^ix_1^j\quad & j \ge 0,\\
x_0^ix_{-1}^{-j}\quad &j <0,
\end{array}\right.\quad i \in \mathbb{N},j \in
\mathbb{Z}  
$$
form a vector space basis of $A$.

We denote by $G$ the 
automorphism group of $A$. 
The defining relations imply that 
for any $\lambda \in \k \!\setminus\!\!\{0\}$
there is a unique  
$\sigma_\lambda \in G$ 
with $\sigma_\lambda (x_n)=\lambda^n x_n$.
\begin{lemma}
Any $\sigma \in G$ is of the form 
$\sigma_\lambda$ for some $ \lambda $,
that is, $G \simeq \k \!\setminus\!\!\{0\}$.
\end{lemma}
\begin{pf}
It is straightforward to classify the characters 
of $A$ and to see that the intersection of their 
kernels is the ideal generated by 
$x_0$. It follows that any automorphism 
$ \sigma $ maps $x_0$ to a nonzero 
scalar multiple of $x_0$. Hence 
$x_0 \sigma(x_{\pm 1})=
q^{\pm 2}\sigma (x_{\pm 1})x_0$ which implies
$\sigma (x_{\pm 1})=f_\pm x_{\pm 1}$
for some $f_\pm \in k[x_0]$.
Inserting into the defining relations 
gives the claim.  
\end{pf}

See 
also \cite{dijkhuizen,chef} for more information
about this algebra. 

\section{The fundamental class} 
As shown in \cite{ulineu}, $A$
satisfies (\ref{podu}) with $ \mathrm{dim}(A)=2$
(as it probably should be for a quantum 2-sphere),
and with $\sm$ determined uniquely by
$$
		  \sm(x_n)=q^{2n} x_n.
$$
So $H_2(A,{}_\sm A) \simeq H^0(A,A)$,
the centre of $A$. This consists only of the 
scalars, hence $[\d A]$ is unique
up to normalisation. 
Hadfield has computed explicit vector space 
bases of all $H_n(A,{}_\sigma A)$ for general 
automorphisms $ \sigma $ \cite{tom} and 
has given
in particular an explicit cycle representing 
$[\d A] \in H_2(A,{}_\sm A)$:
\begin{eqnarray}
		  \d A
&:=& 2  x_1 \otimes 
		  (x_{-1} \otimes x_0-q^2x_0 \otimes x_{-1}) \nonumber\\ 
&& + 2 x_{-1} \otimes (q^{-2} 
		  x_0 \otimes x_1-x_1 \otimes x_0)\nonumber\\ 
&&+1 \otimes (qx_1 \otimes x_{-1}-q^{-1}x_{-1} \otimes x_1
		  +(q-q^{-1})x_0 \otimes x_0) \nonumber\\ 
&&+2x_0 \otimes (x_1 \otimes x_{-1}-x_{-1} \otimes x_1
		  +(q^2-q^{-2})x_0 \otimes x_0).
		  \nonumber 
\end{eqnarray}

\section{The twisted centre and 
$H_2(A,{}_\sigma A)$} 
As we pointed out above, 
the fundamental class is
classically (meaning for the 
coordinate ring of a
smooth variety) represented by a 
nowhere vanishing
algebraic differential form 
of top degree, and this
form can be multiplied by any 
regular function to give a new 
top degree form. Analogously 
we can act via the
cap product by any 
twisted central element 
$a \in H^0(A,{}_\tau A)$ on the
fundamental class $[\d A] \in H_2(A,{}_\sm A)$ 
and obtain another homology 
class in 
$H_2(A,{}_{\tau \circ \sm} A)$. This operation
clarifies completely the 
structure of all the other
nonvanishing 
$H_2(A,{}_\sigma A)$ that were computed by
Hadfield, since they become altogether 
identified with a free $k[x_0]$-module of rank 1:
\begin{lemma}
The twisted centre $\Lambda^0_G(A)$ 
of $A$ is the subalgebra
generated by 
$x_0 \in H^0(A,{}_\sm A)$, and 
for every $[\omega] \in H_2(A,{}_\sigma A)$, 
$\sigma$ any automorphism, there exists 
exactly one polynomial $f \in k[x_0]$ such that 
$[\omega]=[\d A] \smallfrown f$. 
\end{lemma}
\begin{pf}
Since any automorphism fixes $x_0$, 
any twisted central
element must commute with $x_0$ and is hence a
polynomial in $x_0$ (use the vector space basis
$e_{ij}$). Conversely, it is clear that 
$x_0 \in H^0(A,{}_\sm A)$ and hence 
$x_0^i=x_0 \smallsmile \ldots \smallsmile x_0 
\in H^0(A,{}_{\sigma_{q^{2i}}} A)$ (recall Lemma~\ref{relands}). 
The second part follows by Poincar\'e duality (\ref{podu}),
but of course also from Hadfield's explicit
 computations of all the nontrivial 
$H_2(A,{}_\sigma A)$.
\end{pf}

\section{Twisted traces and 
$H_0(A,{}_\sigma A)$}\label{traces}
As a vector space, 
$H_0(A,{}_\sigma A)$ can be described 
as follows \cite{tom}:
\begin{lemma}\label{hzero}
For $ \sigma = \sigma_\lambda$, 
the following is a vector space basis of 
$H_0(A,{}_\sigma A)$:
$$
		  \{[1]\} \cup 
		  \{[x_{\pm 1}^j]\,|\,j \neq 0,\lambda = 1\} \cup 
		  \{[x_0]\,|\,\lambda \neq q^{2i},i >0\} \cup 
		  \{[x_0^i]\,|\,\lambda = q^{2i},i >0\}. 
$$
\end{lemma}
\begin{pf}
By definition, $H_0(A,{}_\sigma A)$ is as a vector
space the quotient of $A$ by the subspace spanned by 
elements of the form $\b(a,b)=ab-\sigma(b)a$.
Since 
$$
		  a \otimes bc=ab \otimes c+\sigma(c)a \otimes b
		  -\b(a,b,c),
$$
one has $\b(a,bc)=\b(ab,c)+\b(\sigma (c)a,b)$,
so $ \mathrm{im}\, \b$ is spanned by the elements
$b^k_{ij}:=\b(e_{ij},x_k)$, $i \in \mathbb{N}$, $j \in \mathbb{Z}$,
$k=-1,0,1$ which are for $ \sigma=\sigma_\lambda$ given by
\begin{eqnarray}
&& b^{-1}_{ij} =
		  (1-\lambda^{-1}q^{2i})e_{ij-1},\quad j \le 0\nonumber\\ 
&& b^{-1}_{ij} =
		  (q^{-4j+2}-\lambda^{-1}q^{2i+2})e_{i+2j-1}
		  +(q^{-2j+1}-\lambda^{-1}q^{2i+1})e_{i+1j-1},\quad j>0,\nonumber\\ 
&& b^0_{ij} = (q^{-2j}-1)e_{i+1j},\nonumber\\  
&& b_{ij}^1 = (1-\lambda q^{-2i})e_{ij+1},\quad j
 \ge 0\nonumber\\
&& b_{ij}^1 = (q^{-4j-2}-\lambda q^{-2i-2})e_{i+2j+1}
+(q^{-2j-1}-\lambda q^{-2i-1})e_{i+1j+1},\quad j<0.\nonumber
\end{eqnarray}
Reducing this by sheer inspection gives that $ \mathrm{im}\, \b$ is
spanned by the elements
\begin{eqnarray}
&& e_{i+1j},\quad (\lambda-1)e_{0j},\quad j \neq 0,\nonumber\\
&& (\lambda-q^{2i+4})q^{-2i-2}e_{i+20}
+(\lambda-q^{2i+2})q^{-2i-1}e_{i+10},\quad i \ge 0.\nonumber
\end{eqnarray}
The claim follows easily.
\end{pf}

Dually, $H_0(A,{}_\sigma A)$ can be described 
in terms of 
$\sigma$-twisted traces, that is, 
linear functionals 
$\int : A \rightarrow k$ satisfying 
$$
		  \int a b = \int \sigma (b)a,\quad
		  a,b \in A.
$$    
Such traces obviously descend to 
well-defined functionals on 
$H_0(A,{}_\sigma A)$ which we denote for
simplicity by the same symbol. 
For each of the basis elements in 
Lemma~\ref{hzero}, we can (and do) 
define one such trace
\begin{eqnarray}
&& \int_{[x_{\pm 1}^j]} e_{kl} :=
		  \left\{
			\begin{array}{ll}
			 1 \quad & k=0,\pm j=l,\\
			 0 \quad & \mbox{otherwise},
			\end{array}\right.\quad j \ge 0,\nonumber\\ 
&& \int_{[x_0]} e_{kl} := \left\{
									\begin{array}{ll}
			1 \quad & k=1,l=0,\\						 
		  (-1)^{k+1}q^{1-k}
		  \frac{1-\lambda q^{-2}}{1-\lambda
		  q^{-2k}}\quad & k>1,l=0,\\
									0 \quad & \mbox{otherwise},
									\end{array}\right.
		  \nonumber\\ 
&& \int_{[x_0^i]} e_{kl} := \left\{
									  \begin{array}{ll}
										1 \quad & k=i,l=0,\\
										0 \quad & \mbox{otherwise},
									  \end{array}\right.\quad i>1.\nonumber  
\end{eqnarray} 
Note that $\int_{[x_0]}$ is defined in such 
a way that the case $\lambda=q^2$ is 
included. Note also that $\int_{[1]}$
is in fact the character $\varepsilon$ 
determined by $ \varepsilon (x_n)=0$. 
Any automorphism of $A$ leaves 
$ \mathrm{ker}\,\varepsilon$ invariant, so this is 
a twisted trace with respect to all 
automorphisms of $A$.
Since we have
for all elements $[\omega],[\eta]$ of the basis of 
$H_0(A,{}_\sigma A)$ from Lemma~\ref{hzero}
$$
		  \int_{[\omega]} [\eta]=\left\{
\begin{array}{ll}
1 \quad & [\omega]=[\eta],\\
0 \quad & \mbox{otherwise},
\end{array}\right.  
$$
we can (and will) 
use the $\int_{[\omega]}$ to determine
the homology class of a given 0-cycle.
For example, it helps describing the cap product action of the
twisted centre on 
$H_0(A,{}_\sigma A)$:
\begin{lemma}
The action of 
$\Lambda^0_G(A)$ on
$H_0(A,{}_{\sigma_\lambda} A)$, is determined by
\begin{eqnarray}
&&  [x_{\pm 1}^j] \smallfrown x_0=0,\quad j
 >0,\nonumber\\
&& [x_0^i] \smallfrown x_0=[x_0^{i+1}]=\left\{
		  \begin{array}{ll}
			0 \quad & i=0,\lambda=q^{2k},k>0,\\
			-q
		  \frac{1-\lambda}{q^2-\lambda}[x_0]\quad & i=1,\lambda \neq q^2.
		  \end{array}\right.\nonumber
\end{eqnarray} 
\end{lemma}
\begin{pf}
As we remarked in (\ref{spaetzle}), $x_0$ acts on
cycles by multiplication from the left. The claim
follows by applying all the above 
constructed twisted
traces to the resulting cycles.  
\end{pf}

In particular, the span of the classes 
$[x_0^i] \in H_0(A,{}_{\sigma^i_\mathrm{mod}}A)$ is the
orbit of $[1] \in H_0(A,A)$ under
the action of the twisted centre. In a sense, the 
span of the $[x_{\pm 1}^j]$ can be viewed as the orbit
of $[1]$ under the cap product with $x_{\pm 1}$, although 
the latter do not belong to the twisted centre of
$A$: for any subalgebra $B \subset A$ with 
$ \sigma (B) \subset B$, there is a map
$H_n(B,{}_\sigma B) \rightarrow H_n(A,{}_\sigma A)$
given on the level of cycles by the embedding of $B$
into $A$ in each tensor component. If 
a class is in the image of this map, then taking the
cap product with a twisted central element of $B$ is
well-defined, and this applies here to the case
$B$ is the subalgebra generated by $x_1$ or $x_{-1}$, respectively.

\section{Three twisted derivations}  
The Podle\'s sphere is a module algebra over 
the Hopf dual $\mathbb{C}_q[SU(2)]^\circ$ of 
the quantised coordinate ring of 
$SU(2)$, hence twisted primitive
elements therein act as twisted derivations on 
$A$. We will not need Hopf algebra 
theory later, so we rather state the following 
lemma that the reader can
verify directly by checking 
compatibility with the
defining relations of $A$: 
\begin{lemma}
The assignments
\begin{eqnarray}
&& \partial_1 : x_{-1},x_0,x_1 \mapsto 
		  0,qx_{-1},1+(q+q^{-1})x_0,\nonumber\\ 
&& \partial_0 : x_{-1},x_0,x_1 \mapsto 
		  -x_{-1},0,x_1,\nonumber\\ 
&& \partial_{-1} : x_{-1},x_0,x_1 \mapsto 
		  1+(q+q^{-1})x_0,q^{-1}x_1,0\nonumber
\end{eqnarray} 
can be extended uniquely to
1-cocycles $ \partial_i \in C^1_{\sigma^{-|i|}_\mathrm{mod}}$,
$$
		  \partial_i(ab)=\sigma^{-|i|}_\mathrm{mod}(a) \partial_i(b)+
		\partial_i (a) b,\quad a,b \in A,\quad i=-1,0,1.
$$
\end{lemma}

The following 
lemma describes the cup product action of
$\Lambda^0_G(A)$ on these derivations. 
Admittedly, the result is slightly weird:
\begin{lemma}
The $\Lambda^0_G(A)$-module generated by 
$[\partial_0]$ is free, but  
$[\partial_{\pm 1}] \smallsmile x_0=0$.
\end{lemma}
\begin{pf}
One checks directly that 
one has for all $a \in A$
$$
		  (\partial_{\pm 1} \smallsmile x_0)(a)=
		  x_0 \partial_{\pm 1}(a)=
		  \pm \frac{1}{q-q^{-1}} (x_{\mp 1}a
		  -ax_{\mp 1}),
$$ 
so the derivations 
$\partial_{\pm 1} \smallsmile x_0$ are
inner. On the other hand, the 
computation of $e_{jk} x_1-
\sigma^i_\mathrm{mod}(x_1)e_{jk}=
b^1_{jk}$  in the proof of Lemma~\ref{hzero}
shows that no inner derivations 
in $C^1_{\sigma^i_\mathrm{mod}}$ can map 
$x_1$ to 
$x_0^ix_1=(\partial_0 \smallsmile x_0^i)(x_1)$.
\end{pf}

For the reason explained at the end of the
previous section, it does make sense to take the cup
product between $ \partial_{\pm 1}$ and  
$x_{\mp 1}^j$, although these are not twisted central, 
and this produces new twisted
derivations. Acting with them on the fundamental class and
comparing the result with the generators of 
$H_1(A,{}_\sigma A)$ computed in \cite{tom}   
allows one to describe all twisted derivations
of $A$, see \cite{tomuli2} where we carried this out
for $ \mathbb{C}_q[SU(2)]$. However, there is little 
gain in this for the main purpose of the
present paper which is to obtain 
a functional describing $\d A$ in a dual fashion, so I
leave out these calculations.

Let us compute instead the algebra generated by the 
$[\partial_i]$.
Classically, a differential 
form can be contracted with a
vector field to reduce its degree, and in the
quantum case this is generalised by 
the cap product action of (the
cohomology class of) a 
twisted derivation on a homology
class. 
Here is the full orbit of 
$\d A$ under this action of
the $ \partial_i$:
\begin{eqnarray}
 \d A \smallfrown \partial_0
&=& 2q^{-2} x_{-1}x_0 \otimes x_1 
		  +2q^2x_1x_0 \otimes x_{-1} 
		  -2(q^2+q^{-2})x_0^2 \otimes x_0 
		  \nonumber\\ 
&&+	  q^{-1}x_{-1}  \otimes x_1 
		  +qx_1  \otimes x_{-1} 
		  -2(q+q^{-1})x_0  \otimes x_0 ,\nonumber\\ 
		  (\d A \smallfrown \partial_0) 
		  \smallfrown \partial_0 
&=& 2(q^2-q^{-2})x_0^3+3(q-q^{-1})x_0^2,
		  \nonumber\\  
		  (\d A \smallfrown \partial_0) 
			\smallfrown \partial_{-1} 
&=& 2(-q^5+q^{-1})x_1x_0^2
		  +(-2q^2+1+q^{-2}1)x_1x_0 +
		  q^{-1}x_1,\nonumber\\
		  (\d A \smallfrown \partial_0) 
			\smallfrown \partial_1
&=& 2(q-q^{-5})x_{-1}x_0^2+
		  (q^2+1-2q^{-2})x_{-1}x_0
		  +qx_{-1},\nonumber
\end{eqnarray}  
\begin{eqnarray}
		  \d A \smallfrown \partial_{-1} 
&=& -2q^{-1}x_1^2 \otimes x_{-1} 
		  -2q^{-1}x_0^2 \otimes x_1 
		  +2(q^3+q^{-3}) x_1 x_0 \otimes x_0 
		  \nonumber\\
&&-(1+q^{-2})x_0 \otimes x_1 
		  +(1+q^{-2})x_1 \otimes x_0 -
		  q^{-1} \otimes x_1,\nonumber\\ 
		  (\d A \smallfrown \partial_{-1}) 
			\smallfrown \partial_0 
&=& 2(q^{-3}-q^3) x_1x_0^2
		  +(-q^2-1+2q^{-2})x_1x_0
		  -q^{-1}x_1,\nonumber\\ 
		  (\d A \smallfrown \partial_{-1}) 
			\smallfrown \partial_{-1} 
&=&2(q^2-q^{-6})x_1^2x_0
		  +(q^{-3}-q^{-5})x_1^2,\nonumber\\ 
		  (\d A \smallfrown \partial_{-1}) 
			\smallfrown \partial_1 
&=& 2(q^{-8}-1)x_0^3
		  +(-q-2q^{-1}+q^{-5}+
		  2q^{-7})x_0^2 \nonumber\\ 
&&+(-2-q^{-2}+q^{-4})x_0 -q^{-1},\nonumber
\end{eqnarray}  
\begin{eqnarray}
		  \d A \smallfrown \partial_1 
&=& 2q x_{-1}^2 \otimes x_1 
		  +2qx_0^2 \otimes x_{-1} 
		  -2(q^3+q^{-3})x_{-1}x_0 \otimes x_0 \nonumber\\ 
&&+(q^2+1)x_0  \otimes x_{-1} 
		  +(-q^2-1)x_{-1}  \otimes x_0 +q \otimes x_{-1},\nonumber\\ 
 (\d A \smallfrown \partial_1) \smallfrown \partial_0 
&=& 2(q^3-q^{-3})x_{-1}x_0^2
		  +(2q^2-1-q^{-2})x_{-1}x_0-qx_{-1},\nonumber\\ 
(\d A \smallfrown \partial_1) \smallfrown \partial_{-1} 
&=& 2(1-q^8)x_0^3
		  +(-2q^7-q^5+2q+q^{-1}) x_0^2 \nonumber\\ 
&& +(-q^4+q^2+2)x_0+q,\nonumber\\ 
(\d A \smallfrown \partial_1) \smallfrown \partial_1 
&=& 2(q^6-q^{-2})x_{-1}^2x_0
		  +(q^5-q^3)x_{-1}^2.\nonumber
\end{eqnarray} 
In homology, this reduces to:
\begin{eqnarray}
&& ([\d A] \smallfrown [\partial_i]) 
			\smallfrown [\partial_i] = 0,\nonumber\\  
&& ([\d A] \smallfrown [\partial_0]) \smallfrown 
		  [\partial_{-1}] = 
		  -([\d A] \smallfrown [\partial_{-1}]) 
			\smallfrown [\partial_0] =
			q^{-1}[x_1],\nonumber\\ 
&& ([\d A] \smallfrown [\partial_0]) 
		  \smallfrown [\partial_1] = 
		  -([\d A] \smallfrown [\partial_1]) 
			\smallfrown [\partial_0] =  
			q[x_{-1}],\nonumber\\ 
&& ([\d A] \smallfrown [\partial_1]) 
			\smallfrown [\partial_{-1}] =
		  -q^2([\d A] \smallfrown [\partial_{-1}]) 
		  \smallfrown [\partial_1] =
		  (q^2+1)[x_0]+q.\nonumber
\end{eqnarray} 
From this, we see:
\begin{lemma}
The $\partial_i$ satisfy no other 
relations than those dictated by 
Lemma~\ref{relands},
$$
		  [\partial_i] \smallsmile [\partial_j] = 
		  -q^{2ij}[\partial_j] \smallsmile 
		  [\partial_i],\quad i \le j. 
$$
\end{lemma}
\begin{pf}
Poincar\'e duality tells that the
$\Lambda_G(A)$-module 
$ \Omega^G(A)$ is free and 
generated by $\d A$, so 
the $ [\partial_i]$ satisfy under $\smallsmile$ all 
the relations 
they satisfy as linear maps on
$\Omega^G(A)$. The result follows 
from the above computations.
\end{pf}
  
\section{The volume form} 
Finally, we now put
$$
		  \varphi : \Omega^G(A) \rightarrow k,
		  \quad
		  [\omega] \mapsto 
		  q^{-1} \int_{[1]} [\omega] \smallfrown 
 		  ([\partial_1] \smallsmile [\partial_{-1}]),
$$
or explicitly on chain level
$$
		  \varphi (a_0,a_1,a_2)=
		  q^{-1} \int_{[1]} 
		  \sigma^{-2}_\mathrm{mod}(a_0)
		  \sigma^{-1}_\mathrm{mod}(\partial_1(a_1))
		  \partial_{-1}(a_2).
$$
Our above computation of 
$([\d A] \smallfrown [\partial_1]) 
\smallfrown [\partial_{-1}]$ implies 
$$
		  \varphi (\d A)=1,
$$
so the functional $\varphi$
provides a dual description of 
the fundamental class. 
It is a very useful  
tool (and probably 
the only really applicable one) for 
checking whether or not a given 
$\sm$-twisted 2-cycle 
has trivial class in $H_2(A,{}_\sm A)$ or 
not, which as a result of the above 
computations (recall that 
$H_2(A,{}_\sm A) \simeq \mathbb{C}$) 
is the fact if and only if 
$ \varphi$ vanishes on the 
cycle.

Note that we remarked in Section~\ref{traces}
that $\int_{[1]}$ is actually a character that 
we also denote by $\varepsilon$ (when 
embedding $A$ as usual into the quantised
coordinate ring of $SU(2)$, 
this character becomes the restriction of 
the counit). Hence $ \varphi $ can also be 
written as simple as
$$
		  \varphi (a_0,a_1,a_2)=
		  q^{-1}\varepsilon(a_0)F(a_1)E(a_2),
$$
where $E,F : A \rightarrow k$ are 
the (untwisted) derivations given by
$$
		  E(a):=\varepsilon(\partial_{-1}(a)),
		  \quad
		  F(a):=\varepsilon(\partial_1(a)).
$$
A moment's reflection now gives 
the explicit formula
\begin{equation}\label{shock}
		  \varphi (e_{ij},e_{kl},e_{mn})=
		  q^{-1} \delta_{i0}\delta_{j0}
		  \delta_{k0}\delta_{l1}\delta_{m0}
		  \delta_{n-1}.
\end{equation} 
This looks surprisingly simple (not to say 
banal), usually one 
expects $\int$ to be 
some sort of integral. However, we
observed already in \cite{tomuli2} 
for the case of quantum 
$SU(2)$ that the functional appearing when 
expressing the volume form dual to 
$\d A$ as above is given by something like 
the integral of the restriction of functions to
a maximal torus which is a Poisson subgroup 
of the Poisson group quantised by 
$\mathbb{C}_q[SU(2)]$. Here we have the 
same phenomenon, 
but it appears much sharper since the functional $\int_{[1]}$ is 
really integration over (meaning evaluation in)
a single point. It somehow 
seems that the homological
information about a quantum space can be 
supported in 
a classical subspace of smaller dimension.

Note also that one can alternatively define
$$
		  \varphi_\pm : \omega \mapsto
		  \pm q^{\mp 1} \int_{[x_{\mp 1}]}
		  \omega \smallfrown (\partial_0
		  \smallsmile \partial_{\pm 1}),
$$
where $\int_{[\pm 1]}$ is the trace (no twist)
dual to $[x_{\pm 1}] \in H_0(A,A)$. Again,
the above computations give
$\varphi_\pm(\d A)=1$,
hence $\varphi=\varphi_+=\varphi_-$ 
as functionals on 
$H_2(A,{}_\sm A) \simeq \mathbb{C}$ 
(but not as 
functionals on $C_2^\sm$), 
and for some purposes 
this representation of the functional 
might be better suited that $\varphi $.

\section{Cyclic homology}
In the classical case $A=k[X]$,
$X$ a smooth variety, exterior derivation 
$\d$ turns the algebraic differential forms 
$\Omega(X)$ into 
a cochain complex that computes 
the algebraic de 
Rham cohomology $H^n(X)$ of $X$. 
In the noncommutative case, 
Connes' cyclic homology provides a subtle 
analogue of $H^n(X)$. The extension 
of cyclic homology to 
the twisted coefficients ${}_\sigma A$
arose in the work of Kustermans, Murphy 
and Tuset on covariant differential 
calculi over quantum groups \cite{kmt},
but can also be viewed as a special case of 
Connes-Moscovici's Hopf-cyclic 
homology \cite{cm}, see e.g.~our article 
\cite{tomuli2} and the references therein 
for more background.\footnote{As Jack 
Shapiro informed me, he will 
discuss in a forthcoming article 
the corresponding version of 
noncommutative de Rham 
theory \cite{jack}.}

The precise 
relation between 
$HC_n(k[X])$ ($ \sigma = \mathrm{id}$)
and $H^n(X)$ is
\begin{equation}\label{derham}
		  HC_n(k[X]) \simeq 
		  \Omega^n(X)/\mathrm{im}\, \d \oplus
		  H^{n-2}(X) \oplus H^{n-4}(X) \oplus \ldots,
\end{equation} 
so all differential $n$-forms (not only the 
closed ones in $ \mathrm{ker}\, \d$!) have 
classes in $HC_n(k[X])$, and this gives  
a map
$$
		  \mathsf{I} : \Omega^n(X) \simeq
		  H_n(k[X],k[X]) \rightarrow
 		  HC_n(k[X]).
$$
Furthermore, $\d$ (applied to 
$\Omega^n(X)/\mathrm{im}\, \d$)
gives a well-defined map
$$
		  \mathsf{B} : 
 		  HC_n(k[X]) \rightarrow 
		  \Omega^{n+1}(X) \simeq
		  H_{n+1}(k[X],k[X])
$$
which kills all the $H^{n-2i}(X)$ 
summands in (\ref{derham}), and finally
there is
$$
		  \mathsf{S} : HC_n(k[X]) \rightarrow 
		  HC_{n-2}(k[X])
$$
that cuts off the first term 
$\Omega^n(X)/\mathrm{im}\, \d$ 
and leaves the rest
untouched (using the obvious embedding 
$H^{n-2}(X) \rightarrow 
\Omega^{n-2}(X)/\mathrm{im}\, \d$ 
in the next summand).

This whole picture carries over to the general
noncommutative case and becomes condensed 
into Connes' $\mathsf{SBI}$-sequence,
see \cite{loday,weibel} for the details. 
The upshot of this is that 
$HC_n(A)$ contains a whole lot of 
ballast, the really interesting part is only 
the image of the natural map 
$\mathsf{I}$ coming from Hochschild 
homology, which equals the kernel of the 
periodicity map $\mathbf{S}$.

One way to define the cyclic 
theory is in terms of the operator
$$
		  \t : C_n^\sigma \rightarrow 
		  C^\sigma_n : 
		  a_0 \otimes \ldots \otimes a_n \mapsto
		  (-1)^n \sigma(a_n) \otimes a_0 \otimes 
		  \ldots \otimes a_{n-1}.
$$
Its coinvariants 
$C^\sigma_n/\mathrm{im}\,(\mathrm{id}-\t)$
form a quotient complex of $(C^\sigma_n,\b)$
whose homology  
is $HC_n^\sigma(A)$ ($k$ should contain
$ \mathbb{Q} $ and $\sigma$ should
be diagonalisable, otherwise the result might 
be not what one wants). 
As a consequence,
a linear functional
$\psi : C_n^\sigma \rightarrow k$
with $ \psi (\mathrm{im}\, \b)=0$ descends
to a functional on 
$HC^\sigma_n(A)$ if it is invariant under $\t$,
$\psi = \psi \circ \t$. Clearly, $\psi$ 
also induces a functional on 
$H_n(A,{}_\sigma A)$, but this might  
vanish even when the one on
$HC^\sigma_n(A)$ doesn't 
(namely when there exists
$\chi : C_{n-1}^\sigma \rightarrow k$
with $\psi = \chi \circ \b$, but no such 
$\chi$ which is cyclic, 
$\chi=\chi \circ \t$). Such a functional 
on $HC_n^\sigma(A)$ that vanishes on 
$H_n(A,{}_\sigma A)$ 
corresponds to the above described 
ballast in cyclic homology, 
it vanishes in the classical case 
on the leading term 
$\Omega^n(X)/\mathrm{im}\, \d$ of 
$HC_n(k[X])$ and is rather 
a functional on some $HC^\sigma_{n-2k}(A)$, 
$k>0$, that is promoted to a functional on 
$HC^\sigma_n(A)$ using the periodicity 
operation $\mathsf{S}$.

\section{The case of the Podls\'s sphere}  
Let now $A$ be again the Podle\'s sphere.
Schm\"udgen and Wagner have constructed 
in \cite{sw} a nontrivial cyclic 2-cocycle 
on $A$ which later was shown by Hadfield 
to be trivial when viewed on 
Hochschild homology. It is now natural to ask 
whether our volume form $ \varphi $ 
constructed above does also 
give rise to a nontrivial functional 
on cyclic homology.

It is easily checked 
that a functional $\psi : C^\sigma_n 
\rightarrow k$ vanishing on $ \mathrm{im}\, \b$
is cyclic if and 
only if $\varphi (1,a_1,\ldots,a_n)=0$ for 
all $a_1,\ldots,a_n \in A$. Using this one sees that
our $ \varphi $ itself is not cyclic since by
(\ref{shock}) we have
$$
		  \varphi (1,x_1,x_{-1})=q^{-1} \neq 0.
$$
However, we can alter $\varphi$ by a 
coboundary to make it cyclic, 
so the problem is only a matter of 
representing
the functional on $H_2(A,{}_\sm A)$ 
properly:
\begin{lemma}
Define $ \eta := \phi \circ \b : 
C_2^\sm \rightarrow k$, where
the linear functional 
$\phi : A \otimes A \rightarrow k$  
vanishes on all $e_{kl} \otimes e_{mn}$
except on the following ones:
$$
		  \phi (1,x_0):=\frac{1}{q^{-2}-1},\quad
		  \phi (1,x_0^2):=\frac{1}{q-q^{-1}},\quad
		  \phi(x_0,x_0):=\frac{1}{2(q-q^{-1})}.
$$
Then $\eta$ induces the trivial functional
on $H_2(A,{}_\sm A)$ and
$\varphi + \eta$ is cyclic.
\end{lemma}
\begin{pf}
By its definition, $ \eta $ vanishes on 
$ \mathrm{im}\, \b$ and defines the trivial 
functional on $H_2(A,{}_\sm A)$.
By (\ref{shock}), the only chain of the form 
$1 \otimes e_{kl} \otimes e_{mn}$
on which $ \varphi $ does not 
vanish is
$1 \otimes x_1 \otimes x_{-1}$, where 
$ \varphi $ has the value $q^{-1}$, 
and one easily checks using
$$
		  \eta (1,a_1,a_2)= 
		  \phi (a_1,a_2)
		  -\phi (1,a_1a_2)
		  +\phi (\sm(a_2),a_1)
$$
that similarly $ \eta (1,e_{kl},e_{mn})$
vanishes except when 
$e_{kl} \otimes e_{mn}=x_1 \otimes x_{-1}$
and then it equals      $-q^{-1}$. The result 
follows.
\end{pf}


\begin{thebibliography}{99}
\bibitem{bz} K.A.~Brown, J.J.~Zhang,  
\emph{Dualising complexes and 
twisted Hochschild (co)homology for 
noetherian Hopf algebras}, arXiv:math.RA/0603732 (2006).

\bibitem{ce} H.~Cartan, S.~Eilenberg, 
\emph{Homological algebra}, 
Princeton University Press (1956).

\bibitem{cm} A.~Connes, H.~Moscovici, 
\emph{Hopf algebras, cyclic cohomology 
and the transverse index theorem},
Comm. Math. Phys.  {\bf 198} no. 1, 199-246 (1998). 

\bibitem{dijkhuizen} M.~Dijkhuizen, 
 T.~Koornwinder, \emph{Quantum homogeneous 
spaces, duality and quantum $2$-spheres},
Geom. Dedicata {\bf 52} no.~3, 291-315 (1994).

\bibitem{farinati} M.~Farinati,
   \emph{Hochschild duality, localization, and smash products},
   J.~Algebra {\bf 284} no.~1
  415-434 (2005).
 
\bibitem{tom} T.~Hadfield, \emph{Twisted cyclic 
homology of all Podle\'s quantum spheres},
		  J.~Geom.~Phys.  
		  {\bf 57} no.~2, 339-351 (2007). 

\bibitem{tomuli} T.~Hadfield, U.~Kr\"ahmer, 
\emph{On the Hochschild homology of quantum SL(N)},
Comptes Rendus Acad. Sci. Paris, Ser. I {\bf 343}, 9-13 (2006)

\bibitem{tomuli2} T.~Hadfield, U.~Kr\"ahmer,
		  \emph{Twisted Homology of Quantum SL(2) - Part II},
     arXiv:0711.4102 (2007).  

\bibitem{hkr}
G. Hochschild, B. Kostant, A. Rosenberg, 
\emph{Differential forms on regular affine algebras}, 
Trans.~AMS {\bf 102}, 383-408 (1962).

\bibitem{chef} 
A.~Klimyk, K.~Schm\"udgen, \emph{Quantum groups and their
representations}, Springer (1997).

\bibitem{uli} U.~Kr\"ahmer, \emph{Poincare duality 
in Hochschild (co)homology}, Proceedings of 
"New techniques in Hopf algebras and 
graded ring theory" (2006). 

\bibitem{ulineu} U.~Kr\"ahmer, \emph{On the Hochschild
		  (co)homology of quantum homogeneous spaces}, 
		  arXiv:0806.0267 (2008).

\bibitem{kmt} J.~Kustermans, G.~Murphy, L.~Tuset, \emph{Differential calculi over quantum groups and twisted cyclic cocycles,} J. Geom. Phys. {\bf 44}, no. 4, 570-594 (2003). 

\bibitem{loday} J.-L.~ Loday,
\emph{Cyclic homology},
Springer-Verlag (1998).

\bibitem{nesttsygan} R.~Nest, B.~Tsygan,
\emph{On the cohomology ring of an algebra},
Progr. Math. {\bf 172},  337-370 (1999). 

\bibitem{Po} P. Podle\'s,
\emph{Quantum spheres},
Lett. Math. Phys. {\bf 14}, 193-202 (1987).

\bibitem{sw} K.~Schm\"udgen, 
E.~Wagner, \emph{Dirac operator and 
a twisted cyclic cocycle on the standard 
Podle\'s quantum sphere}, J.~Reine Angew.~Math.~{\bf
		  574} 219-235 
(2004).

\bibitem{jack} J.~Shapiro,
\emph{Relations between twisted derivations and
twisted cyclic homology}, in preparation.

\bibitem{vdb} M.~van den Bergh, 
 \emph{A relation between Hochschild homology
and cohomology for Gorenstein rings},
Proc. Amer. Math. Soc. {\bf 126}, no. 5, 1345-1348
(1998). Erratum: Proc. Amer. Math. Soc. 
{\bf 130}, no. 9, 2809-2810 (electronic) (2002).

\bibitem{weibel} C.~Weibel,
\emph{An introduction to homological algebra},
Cambridge Univ. Press (1995).
\end{thebibliography}
\end{document}